\documentclass[twoside,12pt]{article}
\usepackage[colorlinks=true]{hyperref}

\usepackage{braket}

\usepackage{amsmath,amssymb,amsthm,graphicx,booktabs,enumerate,mathrsfs,color,subfigure,multicol,mathdots}

\newtheorem{proposition}{Proposition}[section]
\newtheorem{theorem}[proposition]{Theorem}
\newtheorem{lemma}[proposition]{Lemma}
\newtheorem{corollary}[proposition]{Corollary}
\newtheorem{definition}[proposition]{Definition}

\numberwithin{equation}{section}

\title{ Strict Copositivity for a Class of 3rd Order Symmetric Tensors\thanks{This work was supported by the National Natural Science Foundation of P.R. China (Grant No.12171064), by The team project of innovation leading talent in chongqing (No.CQYC20210309536) and by the Foundation of Chongqing Normal university (20XLB009).}} 
\author{ Min Li, Yisheng Song\thanks{Corresponding author: yisheng.song@cqnu.edu.cn}\\
	{\normalsize School of Mathematical Sciences, Chongqing Normal University} \\
	{\normalsize Chongqing, 401331, P.R. China.}\\
	{\normalsize 1753894377@qq.com (Li); yisheng.song@cqnu.edu.cn (Song)}}
\begin{document}
\date{ }
\maketitle

\begin{abstract}
In this article, we mainly give the strictly copositive conditions of a special class of third order three dimensional symmetric tensors. More specifically, by means of the polynomial decomposition method,  the analytic  sufficient and necessary  conditions are established for checking the strict  copositivity of a 3rd order 3-dimensional symmetric tensor with its entries  in $\{-1,0,1\}$.   Several  strict inequalities of cubic ternary homogeneous polynomials are presented by applying these conclusions. Some criteria which ensure the strict  copositivity of a general 3rd order 3-dimensional tensor are obtained


\end{abstract}


\section{Introduction}

\indent A tensor (hypermatrix) is a generalization of vectors and matrices and is easily understood as a multidimensional array.  Clearly, a vector is a first order tensor and a matrix is a second order tensor. A tensor is called {\bf higher order tensor}  if its order is greater than $2$.  As is well known, the structural tensors is class of  tensors with certain structural properties, which had become one of most important topics in tensor researches. Among them, the copositive tensor is an important subclass, which has wide applications in vacuum stability of the Higgs scalar potential  model  \cite{K2012, K2016,CHQ2018,SL2022,S2022,S2023,LS2022}, polynomial optimization \cite{NYZ2018,QCC2018,QL2017,SQ-2015}, hypergraph theory\cite{CHQ2018, NYZ2018, WCWYZ2023, CW2018,CQW2016,WCWYZ2023}, complementarity problems \cite{H-Q2019,HQ2019,QH2019,S-Q2016} and so on,  and includes positive semi-definite tensors \cite{Q2005,Q2006,Q2007}, P-tensors \cite{SQ2015}, semi-positive tensors \cite{SQ2016}, M-tensors \cite{ZQW2014,DQW2014}, Cauchy tensors \cite{CQ2017}, Hilbert tensors \cite{SQ2014,SQ-2016,MS2017}  as its subclasses. 

In 2013, Qi \cite{Q2013} introduced the copositive concept  to  symmetric tensors, and gave some nice properties of such a tensor.  Song-Qi \cite{SQ-2015} gave several necessary and sufficient conditions of copositive tensors by means of  its principal sub-tensors. The copositive conditions of third order two-dimensional symmetric tensors was provided by Schmidt-He$\beta$ \cite{SH1988} in 1988. Recently, Qi-Song-Zhang \cite{QSZ2022}  presented the strictly copositive conditions of third order two-demensional tensors  by means of Sturm theorem. Qi-Song-Zhang \cite{QSZ2023} provided a  copositive condition of  third order three-dimensional symmetric tensors. Liu-Song \cite{LS2022} gave a analytically sufficient condition of the  copositivity  of  third order three-dimensional symmetric tensors. However, the analytic sufficient and necessary condition of the  copositivity is not bubilt for a higher order tensor with its dimension being greater than $2$ even though its order is $3$ or $4$.

The concept of  copositive matrices was  first proposed by Motzkin \cite{M1952} in 1952. Simpson-Spector \cite{SS1983} and  Hadeler \cite{H1983} and Nadler \cite{N1992} and Chang-Sederberg \cite{CS1994} and Andersson-Chang-Elfving \cite{ACE}  respectively proved  the (strict) copositive conditions  of $2\times2$ and $3\times3$ symmetric matrices  by using different methods of argumentation. Li-Feng \cite{LF1993} proved that the copositivity of a $n\times n$ symmetric matrix is equivalent to the copositivity of the difference between the two $(n-1)\times (n-1)$ symmetric matrices  if the matrix has a row whose off-diagonal elements are all non-positive, which  derives a copositive criteria of  symmetric matrices when $n=3,4.$ In  1969, Baston \cite{B1969} and Haynsworth-Hoffman \cite{HH1969} independently provided an analytic method for determining the copositivity of a class of special $n\times n$ matrices with its entries being $1$ or $-1$.  In 1973, Hoffman-Pereira  \cite{HP1973} established the copositive conditions of symmetric matrices with its entries being $\pm1$ or $0$. 
 
 In this paper, we mainly discuss  analytic sufficient and necessary conditions  of the strict copositivity of a special class of third order three-dimensional symmetric tensors with its entries being $\pm1$ or $0$.  With the help of  these conclusions, we gave some criteria of the strict  copositivity of a general 3rd order 3-dimensional tensor. Moreover, we also provide several  strict inequalities of cubic ternary homogeneous polynomials.

\section{Preliminaries}

 In this paper, we use capital letters $A, B, M,$ ..., for matrices and calligraphic letters $\mathcal{A} $, $\mathcal{B}$, $\mathcal{T}$, ..., for tensors. We denote the set of all $m$th order $n$-dimensional real tensors as $T_{m, n}$, where $m$ and $n$ are positive integers, $m, n \ge 2$, unless otherwise stated. For any tensor $\mathcal{A} = (a_{i_1\cdots i_m}) \in T_{m, n}$ , if its entries $a_{i_1\cdots i_m}$ are invariant under any permutation of its indices, then $\mathcal{A}$ is called a \textbf{symmetric tensor}. The set of all $m$th order $n$-dimensional real symmetric tensors is denoted as $S_{m, n}$.
 We write  $\mathbb{R}_{+}^n$ for the non-negative octant in  Euclidean space, i.e.,  $$\mathbb{R}_{+}^n=\{\textbf{\em x}=(x_1,x_2,\cdots,x_n)^\top; x_1\ge0,x_2\ge0,\cdots,x_n\ge0\}.$$

\begin{definition}
	Let $\mathcal{A}  = (a_{i_1\cdots i_m}) \in T_{m, n}$.  Then $\mathcal{A}$ is said to be 
\begin{itemize}
	\item \textbf{copositive} if for all $\textbf{x} = (x_1, \cdots, x_n)^\top \in \mathbb{R}_{+}^n$, 
	$$\mathcal{A} \textbf{x}^m=\sum_{i_1, \cdots, i_m = 1}^n a_{i_1\cdots i_m}x_{i_1}\cdots x_{i_m} \ge 0;$$
	\item  \textbf{strictly copositive} if  for all  $\textbf{x} = (x_1, \cdots, x_n)^\top \in \mathbb{R}_{+}^n$, $\textbf{x} \not = 0$, 
	$$\mathcal{A} \textbf{x}^m > 0.$$
\end{itemize}	
\end{definition}

In 2022, the strictly copositive conditions of a 3rd order 2-demensional tensor was proved by Qi-Song-Zhang \cite{QSZ2022}  by means of Sturm theorem. In 1988, Schmidt-He$\beta$ \cite{SH1988} provided the copositive conditions of such a tensor.
 
\begin{theorem}(\cite{SH1988,QSZ2022}) \label{thm:22}
Suppose that $\mathcal{A} = (a_{ijk}) \in S_{3, 2}$. Then $\mathcal{A}$ is (strictly) copositive if and only if $a_{111} \ge 0\ (> 0)$, $a_{222} \ge 0\ (> 0)$ and $$\begin{cases}
a_{112} \ge 0, a_{122} \ge 0;\\
4a_{111}a_{122}^3 + 4 a_{112}^3a_{222} + a_{111}^2a_{222}^2 - 6a_{111}a_{112}a_{122}a_{222} - 3a_{112}^2
a_{122}^2 \ge 0 \ (> 0).
\end{cases}$$
\end{theorem}

For a tensor $\mathcal{A}= (a_{ijk}) \in S_{3, 2}$, the following conclusion is easy to show by Theorem \ref{thm:22}. For $\textbf{\em x} =(x_1,x_2)^\top,$ we have   $$ \mathcal{A}\textbf{\em x}^3  = a_{111}x_1^3 +3a_{112}x_1^2x_2+ 3a_{122}x_1x_2^2 + a_{222}x_2^3.$$

\begin{lemma}\label{lem:23}
	Let $\mathcal{A} = (a_{ijk})\in S_{3, 2}$  with its entires $$a_{ijk} \in \{1,0,-1\}.$$ Then $\mathcal{A}$ is strictly copositive if and only if $$a_{111} = a_{222} = 1 \ \mbox{ and }\  a_{112}+a_{122}\ge0.$$
\end{lemma}

{\bf Proof.} Necessity (only if).  Suppose $\mathcal{A} $ is strictly copositive. Then for $\textbf{\em x }=  (1, 0)^\top$ or  $(0, 1)^\top$, we have  
$\mathcal{A} \textbf{\em x}^3  =a_{111} >0$ or  $a_{222} >0$,  and hence, $$a_{111} = a_{222} = 1.$$ For $\textbf{\em x }= (1, 1)^\top$, we have
$$
\mathcal{A} \textbf{\em x}^3  = x_1^3 +3a_{112}x_1^2x_2+ 3a_{122}x_1x_2^2 + x_2^3 =3(a_{112}+a_{122})+2>0,
$$
and hence, $$a_{112}+a_{122}>-\dfrac23.$$
This implies that 
$a_{112}+a_{122}\ge0 $ since $a_{ijk} \in \{1,0,-1\}$.

Sufficiency (if).  It is obvious that the condition that $a_{112}+a_{122}\ge0$ is equivalent to $$\begin{cases}
	a_{112}\in \{1,0\} \mbox{ and }a_{122}\in \{1,0\};
	\\  a_{112}a_{122}=-1.
\end{cases}$$  
So, the proof may be divided majorly into $2$ categories.

(1)   Assume $a_{112}\in \{1,0\} \mbox{ and }a_{122}\in \{1,0\}$.  Then  $a_{112}\ge0, a_{122}\ge0$, and so,
$\mathcal{A} $ is strictly copositive by Theorem \ref{thm:22}.

(2) Assuming that  $a_{112}a_{122}=-1$. Then we have
$$4a_{122}^3 + 4 a_{112}^3 + 1 - 6a_{112}a_{122} - 3a_{112}^2a_{122}^2 = 4 - 4 + 1 + 6 - 3 = 4 > 0,$$
and hence, according to Theorem \ref{thm:22}, $\mathcal{A} $ is strictly copositive. \qed\\

The (strict) copositive conditions  of $3\times3$ symmetric matrices were   respectively proved by Simpson-Spector \cite{SS1983} and  Hadeler \cite{H1983} and Nadler \cite{N1992} and Chang-Sederberg \cite{CS1994} and Andersson-Chang-Elfving \cite{ACE}  using different methods of argumentation.

\begin{lemma} (\cite{SS1983,H1983,N1992,CS1994,ACE})\label{lem:24} If $M=(m_{ij})$ is a symmetric $3\times3$ matrix, then $M$ is (strictly) copositive if and only if 
$$\begin{cases}
	m_{11} \geq 0\ (>0), m_{22} \geq 0\ (>0), m_{33} \geq 0\ (>0),\alpha=m_{12}+\sqrt{m_{11}m_{22}}\geq 0\ (>0),\\
	 \beta=m_{13}+\sqrt{m_{11}m_{33}}\geq 0\ (>0),\gamma =m_{23}+\sqrt{m_{22}m_{33}}\geq 0\ (>0),\\
	 \delta=m_{12}\sqrt{m_{33}}+m_{13}\sqrt{m_{22}}+m_{23}\sqrt{m_{11}}+\sqrt{m_{11}m_{22}m_{33}}+\sqrt{2\alpha\beta\gamma}\geq 0\ (>0).
\end{cases}$$
\end{lemma}

\section{Strict Copositivity of Third Order Three-Dimensional Symmetric Tensors}

Let  $\mathcal{A} = (a_{ijk}) \in S_{3, 3}$. Then for $\textbf{\em x} = (x_1, x_2, x_3)^\top$,
\begin{equation*}
  \begin{aligned}
\mathcal{A} \textbf{\em x}^3 =& \sum_{i,j,k = 1}^3 a_{ijk}x_{i}x_{j}x_{k} \\
=& a_{111}x_{1}^3 + a_{222}x_{2}^3 + a_{333}x_{3}^3 + 3a_{112}x_{1}^2x_{2} + 3a_{122}x_{1}x_{2}^2 + 3a_{113}x_{1}^2x_{3}  \\
&  + 3a_{133}x_{1}x_{3}^2+ 3a_{223}x_{2}^2x_{3} + 3a_{233}x_{2}x_{3}^2 + 6a_{123}x_{1}x_{2}x_{3}
  \end{aligned}
\end{equation*}
\begin{theorem} \label{thm:31} Let  $\mathcal{A}  = (a_{ijk}) \in S_{3, 3}$ with its entries $a_{ijk}\in\{-1,0,1\}$. 
Suppose that there exist $r,s,t\in\{1,2,3\}$ with $r\ne s$, $s\ne t$, $r\ne t$ such that $$ a_{rss} = a_{rtt} =-a_{123} = 1.$$Then $\mathcal{A} $ is  strictly copositive if and only if $$a_{iii} =1, a_{iij}+a_{ijj}\ge0, \mbox{ for all }i,j\in\{1,2,3\}\mbox{ with }i\ne j,$$ and $$\begin{cases}
	a_{rrs}+a_{rrt}\ge0;\\
	a_{stt}=a_{sst}=1, \begin{cases}a_{rrs}=0, a_{rrt}=-1 ;\\
    a_{rrs}=-1, a_{rrt}=0.
\end{cases}
\end{cases} $$
\end{theorem}
{\bf Proof.} Without loss the generality,  we might take $r=1,s=2,t=3$, and hence, $$ a_{122} = a_{133} =1, a_{123} = -1.$$  

Necessity (only if).  Since each principal subtensor of $\mathcal{A}$ is strictly copositive,  it follows from Lemma \ref{lem:23} that the first condition, ``$a_{iii} =1 $ and $a_{iij}+a_{ijj}\ge0, \mbox{ for all }i,j\in\{1,2,3\}\mbox{ with }i\ne j$" is necessary. Next assume $\mathcal{A} $ is strictly copositive, but the remaining condition doesn't hold.  Which is divided into $5$ cases.

\textbf{Case 1. } $a_{112} = a_{113} = -1$.   For $\textbf{\em x} = (2, 1, 1)^\top$, we have
$$
  \begin{aligned}
\mathcal{A} \textbf{\em x}^3 = &x_{1}^3 + x_{2}^3 + x_{3}^3 - 3x_{1}^2x_{2} + 3x_{1}x_{2}^2 - 3x_{1}^2x_{3} + 3x_{1}x_{3}^2 \\&+3a_{223}x_{2}^2x_{3} +3a_{233}x_{2}x_{3}^2 - 6x_{1}x_{2}x_{3}\\
= &8 + 1 + 1-12 + 6-12 +  6 + 3(a_{223}+a_{233})-12\\ =& 3(a_{223}+a_{233})-14 < 0;
  \end{aligned}
$$

\textbf{Case 2. } $a_{112} =0$, $a_{113} = -1$  and there is at least one $0$ in $\{a_{233}, a_{223}\}$. Then $\max\{a_{233},a_{223}\}\le1$ since $a_{233}+a_{223}\ge0$.  For $\textbf{\em x} = (3, 1, 1.5)^\top$, we have
$$
\begin{aligned}
	\mathcal{A} \textbf{\em x}^3 = &x_{1}^3 + x_{2}^3 + x_{3}^3 + 3x_{1}x_{2}^2 - 3x_{1}^2x_{3} + 3x_{1}x_{3}^2\\& +3a_{223}x_{2}^2x_{3} +3a_{233}x_{2}x_{3}^2 - 6x_{1}x_{2}x_{3}\\
	\leq&27 + 1 + 1.5^3+ 9-27\times1.5 +  9\times1.5^2 + 3\times1.5^2-18\times1.5\\=&-0.125 < 0;
\end{aligned}
$$

\textbf{Case 3. } $a_{112} =-1$, $a_{113} = 0$  and there is at least one $0$ in $\{a_{233}, a_{223}\}$. Then $\max\{a_{233},a_{223}\}\le1$ since $a_{233}+a_{223}\ge0$.  For $\textbf{\em x} = (3, 1.5, 1)^\top$, we have
$$
\begin{aligned}
	\mathcal{A} \textbf{\em x}^3 = &x_{1}^3 + x_{2}^3 + x_{3}^3 + 3x_{1}x_{2}^2 - 3x_{1}^2x_2 + 3x_{1}x_{3}^2\\& +3a_{223}x_{2}^2x_{3} +3a_{233}x_{2}x_{3}^2 - 6x_{1}x_{2}x_{3}\\
	\leq&27 + 1.5^3+ 1 + 9\times1.5^2-27 \times1.5+  9 + 3\times1.5^2-18\times1.5\\=&-0.125 < 0;
\end{aligned}
$$

\textbf{Case 4. } $a_{112} =-1$, $a_{113} = 0$  and  $a_{233}a_{223}=-1$. Then $a_{233}+a_{223}=0$.  For $\textbf{\em x} = (2, 1, 1)^\top$, we have
$$
\begin{aligned}
	\mathcal{A} \textbf{\em x}^3 = &x_{1}^3 + x_{2}^3 + x_{3}^3 +3x_{1}x_{2}^2 - 3x_{1}^2x_2 + 3x_{1}x_{3}^2 \\&+3a_{223}x_{2}^2x_{3} +3a_{233}x_{2}x_{3}^2 - 6x_{1}x_{2}x_{3}\\
	\leq&8 + 1+ 1 + 6-12+  6 +3(a_{233}+a_{223})-12=-2 < 0;
\end{aligned}
$$

\textbf{Case 5. } $a_{112} =0$, $a_{113} = -1$  and  $a_{233}a_{223}=-1$.  Similarly, for $\textbf{\em x} = (2, 1, 1)^\top$, we also have
$$
\begin{aligned}
	\mathcal{A} \textbf{\em x}^3 = &x_{1}^3 + x_{2}^3 + x_{3}^3 + 3x_{1}x_{2}^2 - 3x_{1}^2x_3 + 3x_{1}x_{3}^2\\& +3a_{223}x_{2}^2x_{3} +3a_{233}x_{2}x_{3}^2 - 6x_{1}x_{2}x_{3}\\
	\leq&8 + 1+ 1 + 6-12+  6 +3(a_{233}+a_{223})-12=-2 < 0.
\end{aligned}
$$

So, anyway, they contradict with the strict copositivity of  $\mathcal{A}$. 
Hence, the remaining conditions are necessary.

Sufficiency  (If). (1)  Assume $a_{112}+a_{113}\ge0.$  For  $\textbf{\em x} = (x_1 ,x_2 , x_3)^\top \in \mathbb{R}_{+}^3$,  $\mathcal{A} \textbf{\em x}^3$ may be rewritten as follows,
$$
  \begin{aligned}
\mathcal{A} \textbf{\em x}^3 = &a_{111}x_{1}^3 + a_{222}x_{2}^3 + a_{333}x_{3}^3 + 3a_{112}x_{1}^2x_{2} + 3a_{122}x_{1}x_{2}^2 + 3a_{113}x_{1}^2x_{3} \\
& + 3a_{133}x_{1}x_{3}^2 + 3a_{223}x_{2}^2x_{3} + 3a_{233}x_{2}x_{3}^2 + 6a_{123}x_{1}x_{2}x_{3}\\
=&x_{1}(a_{111}x_{1}^2 + 3a_{112}x_{1}x_{2} + 3a_{122}x_{2}^2 + 3a_{113}x_{1}x_{3} + 3a_{133}x_{3}^2 + 6a_{123}x_{2}x_{3})\\
& + a_{222}x_{2}^3  + 3a_{223}x_{2}^2x_{3} + 3a_{233}x_{2}x_{3}^2+ a_{333}x_{3}^3\\
=& x_1 (\textbf{\em x}^\top M \textbf{\em x})+\mathcal{A}' \hat{\textbf{\em x}}^3,
 \end{aligned}
$$
where $\hat{\textbf{\em x}}=(x_2,x_3)^\top$, $M$ is a matrix, \begin{equation}\label{eq:31}
	M=\begin{pmatrix}
		a_{111} & \frac{3}{2}a_{112} & \frac{3}{2}a_{113} \\
		\frac{3}{2}a_{112} & 3a_{122} & 3a_{123} \\
		\frac{3}{2}a_{113} & 3a_{123} & 3a_{133}
	\end{pmatrix},
\end{equation}
and $\mathcal{A}'=(a'_{ijk})\in S_{3,2}$ with its entries \begin{equation}\label{eq:32}a'_{111}=a_{222}, a_{112}'= a_{223}, a_{122}'=a_{233}, a_{222}'=a_{333}.\end{equation}

Obviously, a tensor $\mathcal{A}'$ is strictly copositive by Lemma \ref{lem:23} since $a_{223}+a_{233}\ge0$.  That is,
$$\mathcal{A}' \hat{\textbf{\em x}}^3>0, \mbox{ for all }\hat{\textbf{\em x}}\in\mathbb{R}^2_+\setminus\{0\}.$$
Next, we show the matrix $M$ is copositive by means of Lemma \ref{lem:24}.  In fact, as $a_{111} = a_{222}  = a_{333} = a_{122} = a_{133}  =-a_{123} = 1$ and $a_{112}+a_{113}\ge0$, we have
$$\aligned \alpha=&\frac{3}{2}a_{112} + \sqrt{3a_{111}a_{122}}=\frac{3}{2}a_{112} + \sqrt{3}>0,\\  \beta =&\frac{3}{2}a_{113} + \sqrt{3a_{111}a_{133}}=\frac{3}{2}a_{113} + \sqrt3>0,\\ \gamma =&3a_{123} + 3\sqrt{a_{122}a_{133}}=0,\\
  \delta=&\frac{3}{2}a_{112}\sqrt{3a_{133}} + \frac{3}{2}a_{113}\sqrt{3a_{122}} + 3a_{123}\sqrt{a_{111}} + 3\sqrt{a_{111}a_{122}a_{133}} + \sqrt{2\alpha\beta\gamma}\\
  =&\frac{3\sqrt{3}}{2}(a_{112}+a_{113})\ge0.\endaligned$$
Therefore, it follows from Lemma \ref{lem:24} that  the matrix $M$ is copositive. That is,
$$\textbf{\em x}^\top M \textbf{\em x}\ge0, \mbox{ for all }\textbf{\em x}\in\mathbb{R}^3_+.$$
So, when $x_2=x_3=0$, $\mathcal{A} {\textbf{\em x}}^3=a_{111}x_1^3>0$;
 when $x_2$ and $x_3$ are not simultaneously equal to $0$, $\mathcal{A}' \hat{\textbf{\em x}}^3>0$, and hence $$\mathcal{A} \textbf{\em x}^3 =x_1 (\textbf{\em x}^\top M \textbf{\em x})+\mathcal{A}' \hat{\textbf{\em x}}^3>0.$$
 Thus, $\mathcal{A} $ is strictly copositive.

(2) Assume $a_{233}=a_{223}=1$ and $ \begin{cases}a_{112}=0, a_{113}=-1 ;\\
	a_{112}=-1, a_{113}=0.
\end{cases}$

\textbf{Case 1}. $a_{233}=a_{223}=1$ and $a_{112}=0, a_{113}=-1$.  By this time,  $\mathcal{A} \textbf{\em x}^3$ may be rewritten as follows,
$$
\begin{aligned}
	\mathcal{A} \textbf{\em x}^3 = &x_1^3 + a_{222}x_2^3 + x_3^3 + 3x_{1}x_{2}^2 -3x_{1}^2x_{3} + 3x_{1}x_{3}^2 + 3x_{2}^2x_{3} + 3x_{2}x_{3}^2 -6x_{1}x_{2}x_{3}\\
	=&(x_+x_2+x_3)^3-3x_1^2x_2-6x_1^2x_3-12x_1x_2x_3.
\end{aligned}
$$
Solve the constrained minimum problem:
$$\aligned
\min \ &\ (x_+x_2+x_3)^3-3x_1^2x_2-6x_1^2x_3-12x_1x_2x_3 \\
\mbox{ s. t. } & x_1+x_2+x_3=1, x_1\ge0, x_2\ge0, x_3\ge0 
\endaligned$$
to yield the minimum value $\dfrac1{49}$ at a point $\left(\dfrac47,\dfrac17,\dfrac27\right)$, and  hence, $\mathcal{A}\textbf{\em x}^3>0 \mbox{ for all }\textbf{\em x}\in\mathbb{R}^3_+\setminus\{0\}$. That is, $\mathcal{A}$ is strictly copositive.

\textbf{Case 2}. $a_{233}=a_{223}=1$ and $a_{112}=-1, a_{113}=0$.  The conclusion is easy to obtain using the same proof technique as Case 1, we omit it.

 All in all,  $\mathcal{A}$ is  strictly copositive, which completes the proof.    \qed\\

 When $a_{123} = 0$, the following result may be obtained also.
 
 \begin{theorem} \label{thm:32} Let  $\mathcal{A}  = (a_{ijk}) \in S_{3, 3}$ with its entries $a_{ijk}\in\{-1,0,1\}$. 
 	Suppose that there exist $r,s,t\in\{1,2,3\}$ with $r\ne s$, $s\ne t$, $r\ne t$ such that $$ a_{rss} = a_{rtt} =  1, a_{123} = 0.$$ Then $\mathcal{A} $ is  strictly copositive if and only if  $$a_{iii} =1, a_{iij}+a_{ijj}\ge0, \mbox{ for all }i,j\in\{1,2,3\}\mbox{ with }i\ne j,$$ and $$ a_{rrs}+a_{rrt}\ge-1.$$
 \end{theorem}
 {\bf Proof.}  Let $ a_{122} = a_{133} = 1$ and $ a_{112}+a_{113}\ge-1$ without loss the generality.
 
 Necessity (only if).  Similarly, the necessity of the first condition is obvious by Lemma \ref{lem:23}. Now we only need show that the second condition, $ a_{112}+a_{113}\ge-1$ is necessary. Suppose $\mathcal{A} $ is strictly copositive, but such a condition doesn't hold.  That is,  $ a_{112}=a_{113}=-1.$ Then  for $\textbf{\em x} = (3, 1, 1)^\top$, we have
 $$
 \begin{aligned}
 	\mathcal{A} \textbf{\em x}^3 = &x_1^3 + x_2^3 + x_3^3 - 3x_1^2x_2 + 3x_1x_2^2 - 3x_1^2x_3 + 3x_1x_3^2 +3a_{223}x_2^2x_3 +3a_{233}x_2x_3^2\\
 	= &27 + 1 + 1-27 + 9-27 +  9 + 3(a_{223}+a_{233}) \\
 	=&3(a_{223}+a_{233})-7\leq6-7< 0,
 \end{aligned}
 $$
which contradicts with the strict copositivity of  $\mathcal{A}$.  So, the condition is necessary.
 
 Sufficiency  (If).  For  $\textbf{\em x} = (x_1 ,x_2 , x_3)^\top \in \mathbb{R}_{+}^3$,  rewritting $\mathcal{A} \textbf{\em x}^3$ as follows,
 $$
  	\mathcal{A} \textbf{\em x}^3 = x_1 (\textbf{\em x}^\top M \textbf{\em x})+\mathcal{A}' \hat{\textbf{\em x}}^3,
 $$
 where $\hat{\textbf{\em x}}=(x_2,x_3)^\top$, $M$ is a matrix given by \eqref{eq:31}, 
 and $\mathcal{A}'$ is a tensor given by \eqref{eq:32}. Similarly to Theorem \ref{thm:31}, we only need show $M$ is  copositive. In fact,  we have
 $$\aligned \alpha=&\frac{3}{2}a_{112} + \sqrt{3a_{111}a_{122}}=\frac{3}{2}a_{112} + \sqrt{3}>0,\\  \beta =&\frac{3}{2}a_{113} + \sqrt{3a_{111}a_{133}}=\frac{3}{2}a_{113} + \sqrt3>0,\\ \gamma =&3a_{123} + 3\sqrt{a_{122}a_{133}}=3>0,\\
 \delta=&\frac{3}{2}a_{112}\sqrt{3a_{133}} + \frac{3}{2}a_{113}\sqrt{3a_{122}} + 3a_{123}\sqrt{a_{111}} + 3\sqrt{a_{111}a_{122}a_{133}} + \sqrt{2\alpha\beta\gamma}\\
 =&\frac{3\sqrt{3}}{2}(a_{112}+a_{113})+3+ \sqrt{6\alpha\beta}\geq 3-\frac{3\sqrt{3}}{2}+ \sqrt{6\alpha\beta}>0.\endaligned$$
 Therefore, it follows from Lemma \ref{lem:24} that  the matrix $M$ is strictly copositive, and hence,
 $\mathcal{A} $ is strictly copositive.
   \qed
 
 When $a_{123} = 1$, we have the following conclusion also.
 
  \begin{theorem} \label{thm:33} Let  $\mathcal{A}  = (a_{ijk}) \in S_{3, 3}$ with its entries $a_{ijk}\in\{-1,0,1\}$. 
 	Suppose that there are $i,j,k\in\{1,2,3\}$ with $i\ne j$, $j\ne k$, $i\ne k$, $$ a_{ijj} = a_{ikk} = a_{123} = 1.$$ Then $\mathcal{A} $ is  strictly copositive if and only if $$a_{iii} =1, a_{iij}+a_{ijj}\ge0, \mbox{ for all }i,j\in\{1,2,3\}\mbox{ with }i\ne j.$$
 \end{theorem}
 {\bf Proof.}  The necessity directly follows from Lemma \ref{lem:23}.  Now we show the sufficiency. Without loss the generality, let $a_{122} = a_{133} = 1.$ Similarly, we also have
 $$
 \mathcal{A} \textbf{\em x}^3 = x_1 (\textbf{\em x}^\top M \textbf{\em x})+\mathcal{A}' \hat{\textbf{\em x}}^3,
 $$
 where $\hat{\textbf{\em x}}=(x_2,x_3)^\top$, $M$ is a matrix given by \eqref{eq:31}, 
 and $\mathcal{A}'$ is a tensor given by \eqref{eq:32}.  By this time, 
 $$\aligned \alpha=&\frac{3}{2}a_{112} + \sqrt{3a_{111}a_{122}}=\frac{3}{2}a_{112} + \sqrt{3}\ge\sqrt{3}-\frac{3}{2}>0,\\  \beta =&\frac{3}{2}a_{113} + \sqrt{3a_{111}a_{133}}=\frac{3}{2}a_{113} + \sqrt3\ge\sqrt{3}-\frac{3}{2}>0,\\ \gamma =&3a_{123} + 3\sqrt{a_{122}a_{133}}=6>0,\\
 \delta=&\frac{3}{2}a_{112}\sqrt{3a_{133}} + \frac{3}{2}a_{113}\sqrt{3a_{122}} + 3a_{123}\sqrt{a_{111}} + 3\sqrt{a_{111}a_{122}a_{133}} + \sqrt{2\alpha\beta\gamma}\\
 =&\frac{3\sqrt{3}}{2}(a_{112}+a_{113})+6+ \sqrt{12\alpha\beta}\geq 6-3\sqrt{3}+ \left(\sqrt{3}-\frac{3}{2}\right)\sqrt{12}>0.\endaligned$$
 Therefore, it follows from Lemma \ref{lem:24} that  the matrix $M$ is strictly copositive, and hence,
 $\mathcal{A} $ is strictly copositive.
 \qed\\
 
 When $|a_{ijk}| = 1$,  the following corollary may be established easily.
\begin{corollary} \label{cor:34}
Let $\mathcal{A} = (a_{ijk})$ be a third order three-dimensional symmetric tensor with its entires $|a_{ijk}| = 1$. Suppose there are $r,s,t\in\{1,2,3\}$ with $r\ne s$, $s\ne t$, $r\ne t$ such that $$ a_{rss} = a_{rtt} =  1.$$ Then $\mathcal{A}$ is strictly copositive if and only if $$a_{iii} =1, a_{iij}+a_{ijj}\ge0, \mbox{ for all }i,j\in\{1,2,3\}\mbox{ with }i\ne j,$$ and $$\begin{cases}
a_{123} =  1;\\
a_{123} = - 1\mbox{ and } \begin{cases}a_{rrs} = a_{rrt} =1;\\a_{rrs}  a_{rrt} =-1.\end{cases}
\end{cases}$$
\end{corollary}

Applying Theorems \ref{thm:31}, \ref{thm:32} and \ref{thm:33},  it is easy to obtain the following  ternary cubic  inequalities.
\begin{corollary} \label{cor:35} Suppose $(x_1,x_2,x_3)\ne(0,0,0)$ and $x_i\geq0,\ i=1,2,3$. Then 
$$\aligned (x_1+ x_2 + x_3)^3 >&12x_{1}x_{2}x_{3}+6x_{1}^2x_{2} +6x_{2}^2x_{3};\\
 (x_1+ x_2 + x_3)^3 >&12x_{1}x_{2}x_{3}+6x_{1}^2x_3 +6x_{2}^2x_{3};\\
 (x_1+ x_2 + x_3)^3 >&12x_{1}x_{2}x_{3}+6x_1^2x_3 +3x_1^2x_2;\\
 (x_1+ x_2 + x_3)^3 >&6x_{1}x_{2}x_{3}+3x_{1}^2x_3+6x_1^2x_3 +6x_{2}^2x_{3}; \\
 (x_1+ x_2 + x_3)^3 >&6x_{1}^2x_3+6x_1^2x_3 +6x_{2}^2x_{3}. \endaligned$$
 Furthermore, these strict inequalities still hold if $x_2^2x_3$ and $x_2x_3^2$ interchanges,  or  $x_1$ and $x_2$ and $x_3$  are triadic exchangeable.
\end{corollary}

For a tensor $\mathcal{A} = (a_{ijk}) \in S_{3, 3}$ with positive diagonal entries, we also obtain the following conclusions.
\begin{corollary}\label{cor:36}
	Let  $\mathcal{A} = (a_{ijk}) \in S_{3, 3}$ with $a_{iii}>0$ for all $i\in\{1,2,3\}$.  Assume that there exist $r,s,t\in\{1,2,3\}$ with $r\ne s$, $s\ne t$, $r\ne t$ such that
	$$\aligned
	a_{rss}\ge&\sqrt[3]{a_{rrr}a_{sss}^{2}}, a_{rtt}\ge\sqrt[3]{a_{rrr}a_{ttt}^{2}},\\
	a_{rrs}\ge&-\sqrt[3]{a_{rrr}^2a_{sss}}, a_{rrt}\ge-\sqrt[3]{a_{rrr}^2a_{ttt}},\\
	a_{sst}\ge&-\sqrt[3]{a_{sss}^2a_{ttt}}, a_{stt}\ge\sqrt[3]{a_{sss}a_{ttt}^2},\\
	a_{123}\ge&\sqrt[3]{a_{111}a_{222}a_{333}}.
	\endaligned$$
	Then $\mathcal{A} $ is strictly copositive.
\end{corollary}
{\bf Proof.} Without loss the generality, we might take $r=1,s=2,t=3$. Then for $\textbf{\em x} = (x_1, x_2, x_3)^\top$, we have
$$
	\begin{aligned}
		\mathcal{A} \textbf{\em x}^3 =&a_{111}x_{1}^3 + a_{222}x_{2}^3 + a_{333}x_{3}^3 + 3a_{112}x_{1}^2x_{2} + 3a_{122}x_{1}x_{2}^2 + 3a_{113}x_{1}^2x_{3}    \\
		&+ 3a_{133}x_{1}x_{3}^2+ 3a_{223}x_{2}^2x_{3} + 3a_{233}x_{2}x_{3}^2 + 6a_{123}x_{1}x_{2}x_{3}\\
		=&y_1^3+y_2^3+y_3^3+3b_{112}y_1^2y_2+3b_{122}y_1y_2^2+3b_{113}y_1^2y_3+3b_{133}y_1y_3^2\\&+3b_{233}y_2y_3^2+3b_{223}y_1^2y_3+6b_{123}y_1y_2y_3\\
		\ge&y_1^3+y_2^3+y_3^3-3y_1^2y_2+3y_1y_2^2-3y_1^2y_3+3y_1y_3^2\\
		&+3y_2y_3^2-3y_2^2y_3+6y_1y_2y_3=\mathcal{A}' \textbf{\em y}^3,
	\end{aligned}
$$where
$$\aligned y_1=&\sqrt[3]{a_{111}}x_1,  y_2=\sqrt[3]{a_{222}}x_2, y_3=\sqrt[3]{a_{333}}x_3, \\
b_{112}=&\dfrac{a_{112}}{\sqrt[3]{a_{111}^{2}a_{222}}}\ge-1, b_{122}=\dfrac{a_{122}}{\sqrt[3]{a_{111}a_{222}^{2}}}\ge1,\\ b_{113}=&\dfrac{a_{113}}{\sqrt[3]{a_{111}^{2}a_{333}}}\ge-1, b_{133}=\dfrac{a_{133}}{\sqrt[3]{a_{111}a_{333}^{2}}}\ge1,\\
b_{223}=&\dfrac{a_{223}}{\sqrt[3]{a_{222}^{2}a_{333}}}\ge-1, b_{233}=\dfrac{a_{233}}{\sqrt[3]{a_{222}a_{333}^{2}}}\ge1,\\
b_{123}=&\dfrac{a_{123}}{\sqrt[3]{a_{111}a_{222}a_{333}}}\ge1.
\endaligned$$
It follows from Corollary \ref{cor:34} that $\mathcal{A}' \textbf{\em y}^3>0$ for all $\textbf{\em y}\in\mathbb{R}^3_+\setminus\{0\},$ and hence, $\mathcal{A} \textbf{\em x}^3>0$ for all $\textbf{\em x}\in\mathbb{R}^3_+\setminus\{0\}.$ \qed\\

Similarly, we also have the following corollary by Theorem \ref{thm:33}.
\begin{corollary}\label{cor:37}
	Let  $\mathcal{A} = (a_{ijk}) \in S_{3, 3}$ with $a_{iii}>0$ for all $i\in\{1,2,3\}$.  Assume that $a_{123}\ge\sqrt[3]{a_{111}a_{222}a_{333}}$ and there exist $r,s,t\in\{1,2,3\}$ with $r\ne s$, $s\ne t$, $r\ne t$ such that
	$$\aligned
	a_{rss}\ge&\sqrt[3]{a_{rrr}a_{sss}^{2}}, a_{rtt}\ge\sqrt[3]{a_{rrr}a_{ttt}^{2}},\\
	a_{rrs}\ge&-\sqrt[3]{a_{rrr}^2a_{sss}}, a_{rrt}\ge-\sqrt[3]{a_{rrr}^2a_{ttt}},
	a_{sst}\ge0, a_{stt}\ge0.
	\endaligned$$
	Then $\mathcal{A} $ is strictly copositive.
\end{corollary}
By means of Theorem \ref{thm:31}, the following coclusions are easily showed using the similar proof technique as Corollary \ref{cor:36}.
\begin{corollary}\label{cor:38}
	Let  $\mathcal{A} = (a_{ijk}) \in S_{3, 3}$ with $a_{iii}>0$ for all $i\in\{1,2,3\}$.  Suppose $a_{123}\ge-\sqrt[3]{a_{111}a_{222}a_{333}}$ and that there exist $r,s,t\in\{1,2,3\}$ with $r\ne s$, $s\ne t$, $r\ne t$ such that $$a_{rss}\ge\sqrt[3]{a_{rrr}a_{sss}^{2}}, a_{rtt}\ge\sqrt[3]{a_{rrr}a_{ttt}^{2}}$$ and one of the following five conditions holds,
	\begin{itemize}
		\item  $
		a_{rrs}\ge\sqrt[3]{a_{rrr}^2a_{sss}}, a_{rrt}\ge-\sqrt[3]{a_{rrr}^2a_{ttt}},
		a_{sst}\ge-\sqrt[3]{a_{sss}^2a_{ttt}}, a_{stt}\ge\sqrt[3]{a_{sss}a_{ttt}^2};$
		\item $a_{rrs}\ge -\sqrt[3]{a_{rrr}^2a_{sss}}, a_{rrt}\ge0, 
		a_{sst}\ge\sqrt[3]{a_{sss}^2a_{ttt}}, a_{stt}\ge-\sqrt[3]{a_{sss}a_{ttt}^2};$ 
		\item $a_{rrs}\ge0, a_{rrt}\ge0,
		a_{sst}\ge\sqrt[3]{a_{sss}^2a_{ttt}}, a_{stt}\ge-\sqrt[3]{a_{sss}a_{ttt}^2};$ 
		\item $a_{rrs}\ge\sqrt[3]{a_{rrr}^2a_{sss}}, a_{rrt}\ge-\sqrt[3]{a_{rrr}^2a_{ttt}},
		a_{sst}\ge0, a_{stt}\ge0;$
		\item $a_{rrs}\ge 0, a_{rrt}\ge0, a_{sst}\ge0, a_{stt}\ge0.$
	\end{itemize}
	Then $\mathcal{A} $ is strictly copositive.
\end{corollary}
From Theorem \ref{thm:32}, it is easy to bublid the following conclusion.
\begin{corollary}\label{cor:39}
	Let  $\mathcal{A} = (a_{ijk}) \in S_{3, 3}$ with $a_{iii}>0$ for all $i\in\{1,2,3\}$.  Suppose $a_{123}\ge0$ and that there exist $r,s,t\in\{1,2,3\}$ with $r\ne s$, $s\ne t$, $r\ne t$ such that $$a_{rss}\ge\sqrt[3]{a_{rrr}a_{sss}^{2}}, a_{rtt}\ge\sqrt[3]{a_{rrr}a_{ttt}^{2}}$$ and one of the following two conditions holds,
	\begin{itemize}
		\item $a_{rrs}\ge0, a_{rrt}\ge-\sqrt[3]{a_{rrr}^2a_{ttt}},
		a_{sst}\ge0, a_{stt}\ge0;$
		\item $a_{rrs}\ge-\sqrt[3]{a_{rrr}^2a_{sss}}, a_{rrt}\ge0,
		a_{sst}\ge-\sqrt[3]{a_{sss}^2a_{ttt}}, a_{stt}\ge\sqrt[3]{a_{sss}a_{ttt}^2}.$
	\end{itemize}
	Then $\mathcal{A} $ is strictly copositive.
\end{corollary}

\end{document}